\documentclass[a4paper,11pt]{amsart}
\usepackage{graphicx}
\usepackage{mathptmx}
\usepackage{amsmath}
\usepackage{amssymb}
\usepackage{enumitem}
\usepackage{xcolor}
\usepackage{xparse}
\NewDocumentCommand{\eulerian}{omm}
 {%
  \genfrac<>{0pt}{}{#2}{#3}%
  \IfValueT{#1}{_{\!#1}}%
 }

 \newmuskip\pFqmuskip

\newcommand*\pFq[6][8]{%
  \begingroup 
  \pFqmuskip=#1mu\relax
  \mathchardef\normalcomma=\mathcode`,
  \mathcode`\,=\string"8000
  \begingroup\lccode`\~=`\,
  \lowercase{\endgroup\let~}\pFqcomma
  {}_{#2}F_{#3}{\left(\genfrac..{0pt}{}{#4}{#5}\bigg|#6\right)}%
  \endgroup
}
\newcommand{\pFqcomma}{{\normalcomma}\mskip\pFqmuskip}

\newtheorem{theorem}{Theorem}
\newtheorem{lemma}[theorem]{Lemma}

\begin{document}

\title[Fully degenerate Dowling and fully degenerate Bell polynomials]{Some identities of fully degenerate Dowling and fully degenerate Bell polynomials arising from $\lambda$-umbral calculus}

\author{Yuankui Ma}
\address{School of Science, Xi’an Technological University, Xi’an, 710021, Shaanxi, P. R. China}
\email{mayuankui@xatu.edu.cn}

\author{Taekyun  Kim}
\address{School of Science, Xi’an Technological University, Xi’an, 710021, Shaanxi, P. R. China\\
Department of Mathematics, Kwangwoon University, Seoul 139-701, Republic of Korea}
\email{tkkim@kw.ac.kr}

\author{Hyunseok Lee}
\address{Department of Mathematics, Kwangwoon University, Seoul 139-701, Republic of Korea}
\email{luciasconstant@kw.ac.kr}

\author{Dae San Kim}
\address{Department of Mathematics, Sogang University, Seoul 121-742, Republic of Korea}
\email{dskim@sogang.ac.kr}

\subjclass[2010]{11B73; 11B83; 05A19; 05A40}
\keywords{degenerate Whitney numbers of the second kind; fully degenerate Bell polynomials; fully degenerate Dowling polynomials}

\maketitle

\begin{abstract}
Recently, Kim-Kim introduced the $\lambda$-umbral calculus, in which the $\lambda$-Sheffer sequences occupy the central position. In this paper, we introduce the fully degenerate Bell and the fully degenerate Dowling polynomials, and investigate some properties and identities relating to those polynomials with the help of $\lambda$-umbral calculus. Here we note that the fully degenerate Bell poynomials and the fully degenerate Dowling polynomials are respectively degenerate versions of the Bell polynomials and the Dowling polynomials, of which the latters are the natural extension of the Whitney numbers of the second kind.
\end{abstract}

\section{Introduction}

\indent A finite lattice $L$ is geometric if it is a finite semimodular lattice which is also atomic. For details on lattices, one refers to [17]. Dowling constructed an important finite geometric lattice $Q_{n}(G)$ of rank $n$ (see [3,Theorem 3]) which satisfies the chain condition (see [3 ,Theorem 1]), out of a finite set of $n$ elements and a finite group $G$ of order $m$. This is called Dowling lattice of rank $n$ over a finite group of order $m$. We let the interested reader refer to [3] for the details on the construction of $Q_{n}(G)$ and its many fascinating properties. \\
\indent For a finite geometric lattice $L$ of rank $n$, Dowling [3] defined the Whitney numbers $V_{L}(n,k)$ of the first kind and the Whitney numbers $W_{L}(n,k)$ of the second kind.
In particular, if $L$ is the Dowling lattice $Q_{n}(G)$ of rank $n$ over a finite group $G$ of order $m$, then the Whitney numbers of the first kind $V_{Q_{n}(G)}(n,k)$ and the Whitney numbers of the second kind $W_{Q_{n}(G)}(n,k)$ are respectively denoted by $V_{m}(n,k)$ and $W_{m}(n,k)$. These notations are justified, since the Whitney numbers of both kinds depend only on the order $m$ of $G$.
In Corollary 6.1 of [3], it was shown that, for any fixed group $G$ of order $m$, the Whitney numbers $V_{m}(n,k)$ and $W_{m}(n,k)$ satisfy the following Stirling number-like relations:

\begin{align}
&(mx+1)^{n}=\sum_{k=0}^{n}W_{m}(n,k)m^{k}(x)_{k},\label{1A}\\
&m^{n}(x)_{n}=\sum_{k=0}^{n}V_{m}(n,k)(mx+1)^{k},\quad(\mathrm{see}\ [3]).\label{2A}
\end{align}
where $(x)_{0}=1$,\,\, $(x)_{n}=x(x-1)\cdots(x-n+1), \quad (n \ge 1)$.\\
From \eqref{1A} and \eqref{2A}, it is evident that the Whitney numbers satisfy the orthogonality realtions. \par

Let $r\in\mathbb{N}$. Then, as further generalizations of Whitney numbers of both kinds, the $r$-Whitney numbers of the first kind $V_{m}^{(r)}(n,k)$ and those of the second $W_{m}^{(r)}(n,k)$ are respectively defined by 
\begin{equation}
m^{n}(x)_{n}=\sum_{k=0}^{n}V_{m}^{(r)}(n,k)(mx+r)^{k}, \label{3A}
\end{equation}
and 
\begin{equation}
(mx+r)^{n}=\sum_{k=0}^{n}W_{m}^{(r)}(n,k)m^{k}(x)_{k},\quad(n\ge 0). \label{4A}
\end{equation}
We remark that $r$-Whitney numbers of both kinds and their applications were studied by several authors (see\ [14]). 
Combinatorial interpretations for the $r$-Whitney numbers of both kinds were found by Gyimesi-Nyul in [4]. We state them here for the convenience of the reader. \par

\vspace{0.1in}

Let $0 \le k \le n,\, r \ge 0,\, n + r \ge 1$,\, and let $ m \ge 1.$ Then the $r$-Whitney number of the first kind $V_{m}^{(r)}(n, k)$ is equal to  the number of coloured permutations in
$S_{n+r}$ which are the product of $k + r$ disjoint cycles such that \\
$\bullet$ the distinguished elements $1, \dots, r$ belong to distinct cycles, \\
$\bullet$ the smallest elements of the cycles are not coloured, \\
$\bullet$ an element in a cycle containing a distinguished element is not coloured if there are no smaller numbers on the arc from the distinguished element to this element, \\
$\bullet$ the remaining elements are coloured with $m$ colours. \\
Further, $V_{m}^{(0)}(0,0)$ is defined as $V_{m}^{(0)}(0, 0)=1.$ 

\vspace{0.1in}

Let $0 \le k \le n,\, r \ge 0,\, n + r \ge 1$,\, and let $ m \ge 1.$ Then the $r$-Whitney number of the second kind $W_{m}^{(r)}(n, k)$ is equal to the number of coloured partitions of $\left\{1, \dots, n + r \right\}$ into $k + r$ nonempty subsets such that \\
$\bullet$  the distinguished elements $1, \dots , r$ belong to distinct blocks, \\
$\bullet$  the smallest elements of the blocks are not coloured, \\
$\bullet$  elements in blocks containing a distinguished element are not coloured, \\
$\bullet$  the remaining elements are coloured with $m$ colours. \\
Further,  $W_{m}^{(0)}(0,0)$ is defined as $W_{m}^{(0)}(0, 0) = 1$.

\vspace{0.1in}

\indent Carlitz initiated the study of degenerate versions of some special numbers and polynomials in [1]. It regained the interests of some mathematicians and yielded many interesting results. For some of these, one refers to [5--9, 11 and references therein]. Here we would like to mention only one thing. In [7], we were led to introduce the notions of degenerate Sheffer sequences and $\lambda$-Sheffer sequences, starting from the question that what if we replace the usual exponential function appearing in the generating functions of Sheffer sequences by the degenerate exponential function as in \eqref{7}. \par
The aim of this paper is to introduce the fully degenerate Bell and the fully degenerate Dowling polynomials and investigate some properties and identities for them with the help of $\lambda$-umbral calculus.

In more detail, we find an explicit expression and Dobinski-like formula for the fully degenerate Bell polynomials, and the generating function for the fully degenerate Dowling polynomials. Then, by using the $\lambda$-umbral calculus, we derive a formula expressing any polynomial in terms of the fully degenerate Bell polynomials. Then, by applying this formula to the degenerate Bernoulli polynomials, the $\lambda$-falling factorial polynomials and the degenerate poly-Bell polynomials we represent them in terms of the fully degenerate Bell polynomials. In addition, we express the degenerate Bell polynomials in terms of the degenerate bernoulli polynomials of the second kind.
Then we deduce a formula representing any polynomial in terms of the fully degenerate Dowling polynomials. Applying this formula to the degenerate Bernoulli polynomials and the $\lambda$-falling factorial polynomials, we express them in terms of the fully degenerate Dowling polynomials.
Further, we express the fully degenerate Dowling polynomials in terms of the the fully degenerate Bell polynomials.  \par

\section{Preliminaries}
In this section, we briefly go over the $\lambda$-umbral calculus, for the details of which we let the reader refer to [7]. By the way, [15] is an excellent reference for the usual umbral calculus which corresponds to $0$-umbral calculus.
And then we recall some special polynomials and numbers that are needed in this paper. \par
Let $\mathbb{C}$ be the field of complex numbers. Let
\begin{displaymath}
\mathcal{F}=\bigg\{f(t)=m_{k=0}^{\infty}a_{k}\frac{t^{k}}{k!}\ \bigg|\ a_{k}\in\mathbb{C}\bigg\},
\end{displaymath}
and let
\begin{displaymath}
\mathbb{P}=\mathbb{C}[x]=\bigg\{\sum_{i=0}^{\infty}a_{i}x^{i}\ \bigg|\ \textrm{$a_{i}\in\mathbb{C}$\ with $a_{i}=0$ for all but finite number of $i$}\bigg\}.
\end{displaymath}
Let $\mathbb{P}^{*}$ be the vector space of all linear functionals on $\mathbb{P}$. If $\langle L|p(x)\rangle_{\lambda}$ denotes the action of the $\lambda$-linear functional $L$ on the polynomial $p(x)$, then the vector space operations on $\mathbb{P}^{*}$ are defined by
\begin{displaymath}
\langle L+M|p(x)\rangle_{\lambda}=\langle L|p(x)\rangle_{\lambda}+\langle M|p(x)\rangle_{\lambda},\quad \langle cL|p(x)\rangle_{\lambda}=c\langle L|p(x)\rangle_{\lambda},
\end{displaymath}
where $c$ is a complex constant. \par
For $f(t)\in\mathcal{F}$, with $f(t)=\sum_{k=0}^{\infty}a_{k}\frac{t^{k}}{k!}$, we define the $\lambda$-linear functional on $\mathbb{P}$ by
\begin{equation}
\langle f(t)|(x)_{k,\lambda}\rangle_{\lambda}=a_{k},\quad (k\ge 0),\label{1}
\end{equation}
where
\begin{displaymath}
(x)_{0,\lambda}=1,\quad (x)_{n,\lambda}=x(x-\lambda)\cdots(x-(n-1)\lambda),\quad (n\ge 1).
\end{displaymath}
By \eqref{1}, we get
\begin{equation}
\langle t^{k}|(x)_{n,\lambda}\rangle_{\lambda}=\langle 1|(t^{k})_{\lambda}(x)_{n,\lambda}\rangle_{\lambda}=n!\delta_{n,k},\quad (n,k\ge 0), \label{2}
\end{equation}
where $\delta_{n,k}$ is the Kronecker's symbol. \par
For $k\ge 0$, the $\lambda$-differential operator $t^{k}$ on $\mathbb{P}$ is defined by
\begin{equation}
(t^{k})_{\lambda}(x)_{n,\lambda}=\left\{\begin{array}{ccc}
(n)_{k}x^{n-k}, & \textrm{if $k\le n$,} \\
0, & \textrm{otherwise.}
\end{array}\right.\label{3}
\end{equation}
where $(x)_{0}=1,\ (x)_{n}=x(x-1)\cdots(x-n+1),\ (n\ge 1)$. \par
Let $f(t)=\sum_{k=0}^{\infty}a_{k}\frac{t^{k}}{k!}\in\mathcal{F}$. The order $o(f(t))$ of $f(t)$ is the smallest integer $k$ which $a_{k}$ does not vanish. \par
For $f(t),g(t)\in\mathcal{F}$ with $o(f(t))=1,\ o(g(t))=0$, there exists a unique sequence $s_{n,\lambda}(x)$ of polynomials such that
\begin{equation}
\big\langle 	g(t)(f(t))^{k}\big|S_{n,\lambda}(x)\big\rangle_{\lambda}=n!\delta_{n,k},\quad (n,k\ge 0).\label{4}
\end{equation}
The sequence $S_{n,\lambda}(x)$ is called the $\lambda$-Sheffer sequence for $(g(t),f(t))$, which is denoted by $S_{n,\lambda}(x)\sim (g(t),f(t))_{\lambda}$. \par
For any nonzero $\lambda\in\mathbb{R}$, the degenerate exponential functions are defined as
\begin{equation}
e_{\lambda}^{x}(t)=\sum_{k=0}^{\infty}(x)_{k,\lambda}\frac{t^{k}}{k!},\quad e_{\lambda}(t)=e_{\lambda}^{1}(t),\quad (\mathrm{see}\ [8,11-13]).\label{5}
\end{equation}
Let $\log_{\lambda}(t)$ be the compositional inverse function of $e_{\lambda}(t)$ such that $\log_{\lambda}\big(e_{\lambda}(t)\big)=e_{\lambda}\big(\log_{\lambda}(t)\big)=t$. \par
Then we have
\begin{equation}
\log_{\lambda}(1+t)=\sum_{n=1}^{\infty}\lambda^{n-1}(1)_{n,1/\lambda}\frac{t^{n}}{n!},\quad (\mathrm{see}\ [6]).\label{6}	
\end{equation}
From \eqref{4}, we note that $s_{n,\lambda}(x)\sim (g(t),f(t))_{\lambda}$ if and only if
\begin{equation}
\frac{1}{g(\bar{f}(t))}e_{\lambda}^{x}\big(\bar{f}(t)\big)=\sum_{n=0}^{\infty}s_{n,\lambda}(x)\frac{t^{n}}{n!}, \label{7}
\end{equation}
for all $x\in\mathbb{C}$, where $\bar{f}(t)$ is the compositional inverse function of $f(t)$ such that $f(\bar{f}(t))=\bar{f}(f(t))=t$. \par
For $s_{n,\lambda}(x)\sim (g(t),f(t))_{\lambda}$, $r_{n,\lambda}(x)\sim (h(t),l(t))_{\lambda}$, we have
\begin{equation}
s_{n,\lambda}(x)=\sum_{k=0}^{n}c_{n,k}r_{k,\lambda}(x),\label{8}
\end{equation}
where
\begin{equation}
c_{n,k}=\frac{1}{k!}\bigg\langle\frac{h(\bar{f}(t))}{g(\bar{f}(t))}\big(l(\bar{f}(t))\big)^{k}\ \bigg|\ (x)_{n,\lambda}\bigg\rangle_{\lambda},\quad (n \ge k \ge 0).\label{9}
\end{equation}

\vspace{0.1in}

{\it{Throughout this paper, we assume that $m$ is any fixed positive integer.}} \\
In [9], the degenerate Whitney numbers of the second kind are given by 
\begin{equation}
e_{\lambda}(t)\frac{1}{k!}\bigg(\frac{e_{\lambda}^{m}(t)-1}{m}\bigg)^{k}=\sum_{n=k}^{\infty}W_{m,\lambda}(n,k)\frac{t^{n}}{n!},\quad (k\ge 0).\label{10}
\end{equation}
Recently, Kim-Kim introduced the degenerate Dowling polynomials given by 
\begin{equation}
e_{\lambda}(t)e^{x\big(\frac{e_{\lambda}^{m}(t)-1}{m}\big)}=\sum_{n=0}^{\infty}D_{m,\lambda}(n,x)\frac{t^{n}}{n!},\quad (\mathrm{see}\ [9]). \label{11}
\end{equation}
From \eqref{11}, we note that
\begin{equation}
D_{m,\lambda}(n,x)=\sum_{k=0}^{n}W_{m,\lambda}(n,k)x^{k},\quad (n\ge 0),\quad (\mathrm{see}\ [9]).\label{12}
\end{equation}
In [7], Kim et al. introduced the partial degenerate Bell polynomials which are given by
\begin{equation}
e^{x(e_{\lambda}(t)-1)}=\sum_{n=0}^{\infty}\mathrm{Bel}_{n,\lambda}(x)\frac{t^{n}}{n!}. \label{13}
\end{equation}
Note that $\displaystyle\lim_{\lambda\rightarrow 0}\mathrm{Bel}_{n,\lambda}(x)=\mathrm{Bel}_{n}(x)\displaystyle$. Here $\mathrm{Bel}_{n}(x)$ are the ordinary Bell polynomials given by
\begin{displaymath}
e^{x(e^{t}-1)}=\sum_{n=0}^{\infty}\mathrm{Bel}_{n}(x)\frac{t^{n}}{n!},\quad (\mathrm{see}\ [2,15]).
\end{displaymath}
As is well known, the Stirling numbers of the first kind are defined by
\begin{equation}
(x)_{n}=\sum_{l=0}^{n}S_{1}(n,l)x^{l},\quad(n\ge 0),\quad(\mathrm{see}\ [2,15]).\label{14}
\end{equation}
By the inversion of \eqref{14}, the Stirling numbers of the second kind are defined by
\begin{equation}
x^{n}=\sum_{l=0}^{n}S_{2}(n,l)(x)_{l},\quad(n\ge 0),\quad(\mathrm{see}\ [15,16]). \label{15}
\end{equation}
In [1], Carlitz introduced the degenerate Bernoulli polynomials defined by
\begin{equation}
\frac{t}{e_{\lambda}(t)-1}e_{\lambda}^{x}(t)=\sum_{n=0}^{\infty}\beta_{n,\lambda}(x)\frac{t^{n}}{n!}.\label{16}
\end{equation}
When $x=0$, $\beta_{n,\lambda}=\beta_{n,\lambda}(0)$ are called the degenerate Bernoulli numbers. \par
The degenerate Bernoulli polynomials of the second kind are given by
\begin{equation}
\frac{t}{\log_{\lambda}(1+t)}(1+t)^{x}=\sum_{n=0}^{\infty}b_{n,\lambda}(x)\frac{t^{n}}{n!},\quad (\mathrm{see}\ [10]).\label{17}
\end{equation}
When $x=0$, $b_{n,\lambda}=b_{n,\lambda}(0)$ are called the degenerate Bernoulli numbers of the second kind. \par
For any $k\in\mathbb{Z}$, the degenerate polyexponential function is defined by Kim-Kim as
\begin{equation}
\mathrm{Ei}_{k,\lambda}(x)=\sum_{n=1}^{\infty}\frac{(1)_{n,\lambda}}{(n-1)!n^{k}}x^{n},\quad (\mathrm{see}\ [11]).\label{18}	
\end{equation}
The degenerate poly-Bell polynomials are defined by Kim et al. as
\begin{equation}
\frac{\mathrm{Ei}_{k,\lambda}(\log_{\lambda}(1+t))}{e_{\lambda}(t)-1}e_{\lambda}^{x}(t)=\sum_{n=0}^{\infty}B_{n,\lambda}^{(k)}(x)\frac{t^{n}}{n!},\quad (\mathrm{see}\ [11]).\label{19}	
\end{equation}
When $x=0$, $B_{n,\lambda}^{(k)}=B_{n,\lambda}^{(k)}(0)$ are called the degenerate poly-Bernoulli numbers. \par
The degenerate Stirling numbers of the first kind are defined  by 
\begin{equation}
\frac{1}{k!}\big(\log_{\lambda}(1+t)\big)^{k}=\sum_{n=k}^{\infty}S_{1,\lambda}(n,k)\frac{t^{n}}{n!},\quad (\mathrm{see}\ [6]).\label{20}	
\end{equation}
Note that $\lim_{\lambda\rightarrow 0}S_{1,\lambda}(n,k)=S_{1}(n,k),\ (n \ge k \ge 0)$. \par
By the inversion of \eqref{20}, the degenerate Stirling numbers of the second kind are defined by 
\begin{equation}
\frac{1}{k!}\big(e_{\lambda}(t)-1\big)^{k}=\sum_{n=k}^{\infty}S_{2,\lambda}(n,k)\frac{t^{n}}{n!},\quad(\mathrm{see}\ [6]).\label{21}
\end{equation}
Note that $\displaystyle\lim_{\lambda\rightarrow 0}S_{2,\lambda}(n,k)=S_{2}(n,k)\displaystyle,\ (n \ge k \ge 0)$. \par

\section{Fully degenerate Bell and fully degenerate Dowling polynomials arising from  $\lambda$-umbral calculus}

Here we introduce the fully degenerate Bell and the fully degenerate Dowling polynomials, and investigate some properties and identities relating to them with the help of $\lambda$-umbral calculus. \\
First, we define the {\it{fully degenerate Bell polynomials}} which are given by $\phi_{n,\lambda}(x)\sim(1,\log_{\lambda}(1+t)\big)_{\lambda}$.
Then, by \eqref{7}, we have
\begin{equation}
e_{\lambda}^{x}(e_{\lambda}(t)-1)=\sum_{n=0}^{\infty}\phi_{n,\lambda}(x)\frac{t^{n}}{n!}.\label{23}
\end{equation}
Note that $\displaystyle\lim_{\lambda\rightarrow 0}\phi_{n,\lambda}(x)=\mathrm{Bel}_{n}(x)\displaystyle$. \par
Next, we observe that 
\begin{align}
\sum_{n=0}^{\infty}\phi_{n,\lambda}(x+y)\frac{t^{n}}{n!}&=e_{\lambda}^{x}(e_{\lambda}(t)-1)e_{\lambda}^{y}(e_{\lambda}(t)-1) \label{24}\\
&=\sum_{n=0}^{\infty}\bigg(\sum_{l=0}^{n}\binom{n}{l}\phi_{l,\lambda}(x)\phi_{n-l,\lambda}(y)\bigg)\frac{t^{n}}{n!}. \nonumber	
\end{align}
From \eqref{24}, we note that
\begin{equation}
\phi_{n,\lambda}(x+y)=\sum_{l=0}^{n}\binom{n}{l}\phi_{l,\lambda}(x)\phi_{n-l,\lambda}(y),\quad (n\ge 0).\label{25}
\end{equation}
Further, we see that
\begin{align}
e_{\lambda}^{x}\big(e_{\lambda}(t)-1\big)&=\sum_{k=0}^{\infty}(x)_{k,\lambda}\frac{1}{k!}\big(e_{\lambda}(k)-1\big)^{k} \label{26} \\
&=\sum_{n=0}^{\infty}\bigg(\sum_{k=0}^{n}(x)_{k,\lambda}S_{2,\lambda}(n,k)\bigg)\frac{t^{n}}{n!}.\nonumber
\end{align}
Therefore, by \eqref{23} and \eqref{26}, we obtain the following lemma.
\begin{lemma}
For $n\ge 0$, we have
\begin{displaymath}
\phi_{n,\lambda}(x)=\sum_{k=0}^{n}S_{2,\lambda}(n,k)(x)_{k,\lambda}.	
\end{displaymath}
\end{lemma}
Now, we note that
\begin{align}
e_{\lambda}^{x}\big(e_{\lambda}(t)-1\big)&=(1-\lambda+\lambda e_{\lambda}(t))^{\frac{x}{\lambda}}=(1-\lambda)^{\frac{x}{\lambda}}\bigg(1+\frac{\lambda}{1-\lambda}e_{\lambda}(t)\bigg)^{\frac{x}{\lambda}}\label{27}\\
&=e_{\lambda}^{x}(-1)	\bigg(1+\frac{\lambda}{1-\lambda}e_{\lambda}(t)\bigg)^{\frac{x}{\lambda}}\nonumber \\
&=e_{\lambda}^{x}(-1)\sum_{k=0}^{\infty}(x)_{k,\lambda}\bigg(\frac{1}{1-\lambda}\bigg)^{k}\frac{1}{k!}e_{\lambda}^{k}(t) \nonumber \\
&=\sum_{n=0}^{\infty}\bigg\{e_{\lambda}^{x}(-1)\sum_{k=0}^{\infty}\frac{(k)_{n,\lambda}}{k!}\bigg(\frac{1}{1-\lambda}\bigg)^{k}(x)_{k,\lambda}\bigg\}\frac{t^{n}}{n!}.\nonumber
\end{align}
Therefore, by \eqref{23} and \eqref{27}, we obtain the following theorem.
\begin{theorem}[Dobinski-like formula]
For $n\ge 0$, we have
\begin{displaymath}
\phi_{n,\lambda}(x)= e_{\lambda}^{x}(-1)\sum_{k=0}^{\infty}\frac{(k)_{n,\lambda}}{k!}\bigg(\frac{1}{1-\lambda}\bigg)^{k}(x)_{k,\lambda}.
\end{displaymath}
\end{theorem}
In view of Lemma 1 and \eqref{12}, we define the {\it{fully degenerate Dowling polynomials}} by
\begin{equation}
d_{m,\lambda}(n,x)=\sum_{k=0}^{n}W_{m,\lambda}(n,k)(x)_{k,\lambda},\quad (n\ge 0).\label{28}
\end{equation}
When $x=1$, $d_{m,\lambda}(n,1)=d_{m,\lambda}(n)$ are called the fully degenerate Dowling numbers. \par
From \eqref{10} and \eqref{28}, we note that
\begin{align}
\sum_{n=0}^{\infty}d_{m,\lambda}(n,x)\frac{t^{n}}{n!}&=\sum_{n=0}^{\infty}\bigg(\sum_{k=0}^{n}W_{m,\lambda}(n,k)(x)_{k,\lambda}\bigg)\frac{t^{n}}{n!} \label{29} \\
&=\sum_{k=0}^{\infty}\bigg(\sum_{n=k}^{\infty}W_{m,\lambda}(n,k)\frac{t^{n}}{n!}\bigg)(x)_{k,\lambda}\nonumber \\
&=\sum_{k=0}^{\infty}e_{\lambda}(t)\frac{1}{k!}\bigg(\frac{e_{\lambda}^{m}(t)-1}{m}\bigg)^{k}(x)_{k,\lambda}\nonumber\\
&=e_{\lambda}(t)e_{\lambda}^{x}\bigg(\frac{e_{\lambda}^{m}(t)-1}{m}\bigg). \nonumber
\end{align}
Therefore, by \eqref{29}, we obtain the generating function of fully degenerate Dowling polynomials.
\begin{theorem}
\begin{displaymath}
e_{\lambda}(t)e_{\lambda}^{x} \bigg(\frac{e_{\lambda}^{m}(t)-1}{m}\bigg)=\sum_{n=0}^{\infty}d_{m,\lambda}(n,x)\frac{t^{n}}{n!}.
\end{displaymath}	
\end{theorem}
Note that $\displaystyle \lim_{\lambda\rightarrow 0}d_{m,\lambda}(n,x)=D_{m}(n,x)$. Here $D_{m}(n,x)$ are the ordinary Dowling polynomials given by
\begin{displaymath}
e^{t}e^{x\big(\frac{e^{mt}-1}{m}\big)}=\sum_{n=0}^{\infty}D_{m}(n,x)\frac{t^{n}}{n!}.
\end{displaymath}
Let us assume that
\begin{displaymath}
\mathbb{P}_{n}=\{p(x)\in\mathbb{C}[x]\ |\ \deg p(x)\le n\},\quad (n\ge 0).
\end{displaymath}
Then $\mathbb{P}_{n}$ is an $(n+1)$-dimensional vector space over $\mathbb{C}$. \par
For $p(x)\in\mathbb{P}_{n}$, we let
\begin{equation}
p(x)=\sum_{k=0}^{n}C_{k}\phi_{k,\lambda}(x).\label{30}
\end{equation}
Then, by \eqref{4}, we get
\begin{align}
\big\langle (\log_{\lambda}(1+t))^{k} \big| p(x)\rangle_{\lambda}&=\sum_{l=0}^{n}C_{l}\big\langle (\log_{\lambda}(1+t))^{k}\big|\phi_{l,\lambda}(x)\big\rangle_{\lambda} \label{31} \\
&=\sum_{l=0}^{n}C_{l}k!\delta_{k,l}=k!C_{k}.\nonumber
\end{align}
Thus, by \eqref{31}, we get
\begin{equation}
C_{k}=\frac{1}{k!}\big\langle (\log_{\lambda}(1+t))^{k}\big|p(x)\big\rangle_{\lambda}.\label{32}
\end{equation}
Therefore, by \eqref{30} and \eqref{32}, we obtain the following theorem.
\begin{theorem}
For $p(x)\in\mathbb{P}_{n}$, we have
\begin{displaymath}
	p(x)=\sum_{k=0}^{n}C_{k}\phi_{k,\lambda}(x),
\end{displaymath}	
where
\begin{displaymath}
	C_{k}=\frac{1}{k!}\big\langle (\log_{\lambda}(1+t))^{k}\big|p(x)\big\rangle_{\lambda}.
\end{displaymath}
\end{theorem}
Let us take
\begin{displaymath}
	p(x)=\beta_{n,\lambda}(x)=\sum_{l=0}^{n}\binom{n}{l}\beta_{n-l,\lambda}(x)_{l,\lambda}\in\mathbb{P}_{n}.
\end{displaymath}
Then, by Theorem 4, we get
\begin{equation}
\beta_{n,\lambda}(x)=\sum_{k=0}^{n}C_{k}\phi_{k,\lambda}(x),\label{33}	
\end{equation}
where
\begin{align}
C_{k}&=\frac{1}{k!}\big\langle (\log_{\lambda}(1+t))^{k}\big|\beta_{n,\lambda}(x)\big\rangle_{\lambda}\label{34} \\
&=\frac{1}{k!}\sum_{l=0}^{n}\binom{n}{l}\beta_{n-l,\lambda}\big\langle	(\log_{\lambda}(1+t))^{k}\big|(x)_{l,\lambda}\big\rangle_{\lambda} \nonumber \\
&=\sum_{l=0}^{n}\binom{n}{l}\beta_{n-l,\lambda}\sum_{j=k}^{\infty}S_{1,\lambda}(j,k)\frac{1}{j!}\langle t^{j}|(x)_{l,\lambda}\rangle_{\lambda}\nonumber \\
&= \sum_{l=0}^{n}\binom{n}{l}\beta_{n-l,\lambda}S_{1,\lambda}(l,k).\nonumber
\end{align}
Therefore, by \eqref{33} and \eqref{34}, we obtain the following theorem.
\begin{theorem}
For $n\ge 0$, we have
\begin{displaymath}
\beta_{n,\lambda}(x)=\sum_{k=0}^{n}\bigg\{\sum_{l=0}^{n}\binom{n}{l}\beta_{n-l,\lambda}S_{1,\lambda}(l,k)\bigg\}\phi_{k,\lambda}(x).
\end{displaymath}	
\end{theorem}
Let us take $p(x)=(x)_{n,\lambda}\in\mathbb{P}_{n}$. Then, by Theorem 4, we get
\begin{equation}
(x)_{n,\lambda}=\sum_{k=0}^{n}C_{k}\phi_{k,\lambda}(x),\label{35}	
\end{equation}
where
\begin{align}
C_{k}&=\frac{1}{k!}\big\langle (\log_{\lambda}(1+t))^{k}\big|(x)_{n,\lambda}\big\rangle_{\lambda}	\label{36} \\
&=\sum_{m=k}^{\infty}S_{1,\lambda}(m,k)\frac{1}{m!}\langle t^{m}|(x)_{n,\lambda}\rangle_{\lambda} \nonumber \\
&=S_{1,\lambda}(n,k).\nonumber
\end{align}
Therefore, by \eqref{35} and \eqref{36}, we obtain the following theorem.
\begin{theorem}
For $n\ge 0$, we have
\begin{displaymath}
(x)_{n,\lambda}=\sum_{k=0}^{n}S_{1,\lambda}(n,k)\phi_{k,\lambda}(x).
\end{displaymath}	
\end{theorem}
We note that
\begin{align*}
&B_{n,\lambda}^{(k)}(x)\sim\bigg(\frac{e_{\lambda}(t)-1}{\mathrm{Ei}_{k,\lambda}(\log_{\lambda}(1+t))},t\bigg)_{\lambda}, 	\\
&\phi_{n,\lambda}(x)\sim (1,\log_{\lambda}(1+t))_{\lambda}.
\end{align*}
We assume that
\begin{equation}
B_{n,\lambda}^{(k)}(x)=\sum_{k=0}^{n}C_{n,k}\phi_{k,\lambda}(x). \label{37}	
\end{equation}
Then, by \eqref{8}, we get
\begin{align}
C_{n,k}&=\frac{1}{k!}\bigg\langle\frac{\mathrm{Ei}_{k,\lambda}(\log_{\lambda}(1+t)}{e_{\lambda}(t)-1}\big(\log_{\lambda}(1+t)\big)^{k}\bigg|(x)_{n,\lambda}\bigg\rangle_{\lambda} \label{38} \\
&=\sum_{l=k}^{n}\frac{S_{1,\lambda}(l,k)}{l!} \bigg\langle\frac{\mathrm{Ei}_{k,\lambda}(\log_{\lambda}(1+t)}{e_{\lambda}(t)-1}t^{l}\bigg|(x)_{n,\lambda}\bigg\rangle_{\lambda}\nonumber\\
&=\sum_{l=k}^{n}S_{1,\lambda}(l,k)\binom{n}{l}\sum_{j=0}^{\infty}B_{j,\lambda}^{(k)}\frac{1}{j!}\langle t^{j}|(x)_{n-l,\lambda}\rangle_{\lambda} \nonumber \\
&=\sum_{l=k}^{n}\binom{n}{l}S_{1,\lambda}(l,k)B_{n-l,\lambda}^{(k)}.\nonumber
\end{align}
Therefore, by \eqref{37} and \eqref{38}, we obtain the following theorem.
\begin{theorem}
For $n\ge 0$, we have
\begin{displaymath}
B_{n,\lambda}^{(k)}(x)=\sum_{k=0}^{n}\bigg\{\sum_{l=k}^{n}\binom{n}{l}S_{1,\lambda}(l,k)B_{n-l,\lambda}^{(k)}\bigg\}\phi_{k,\lambda}(x).
\end{displaymath}	
\end{theorem}
From \eqref{17}, we have
\begin{equation}
b_{n,\lambda}(x)\sim\bigg(\frac{t}{e_{\lambda}(t)-1},e_{\lambda}(t)-1\bigg)_{\lambda}.\label{39}
\end{equation}
Thus, by \eqref{8} and \eqref{9}, we get
\begin{equation}
\phi_{n,\lambda}(x)=\sum_{k=0}^{n}C_{n,k}b_{k,\lambda}(x),\label{40}	
\end{equation}
where, for $k=0$,
\begin{align}
C_{n,0}&=\bigg\langle \frac{e_{\lambda}(t)-1}{e_{\lambda}(e_{\lambda}(t)-1)-1}\bigg|(x)_{n,\lambda}\bigg\rangle_{\lambda} \label{41} \\
&=\sum_{l=0}^{n}\frac{\beta_{l,\lambda}}{l!}\big\langle (e_{\lambda}(t)-1)^{l}\big|(x)_{n,\lambda}\big\rangle_{\lambda} \nonumber \\
&=\sum_{l=0}^{n}\beta_{l,\lambda}S_{2,\lambda}(n,l); \nonumber	
\end{align}
for $k\ge 1$,
\begin{align}
C_{n,k}&=\frac{1}{k!}\bigg\langle \frac{e_{\lambda}(t)-1}{e_{\lambda}(e_{\lambda}(t)-1)-1}\big(e_{\lambda}(e_{\lambda}(t)-1)-1\big)^{k}\bigg|(x)_{n,\lambda}\bigg\rangle_{\lambda}	\label{42} \\
&=\frac{1}{k!}\big\langle(e_{\lambda}(t)-1)\big(e_{\lambda}(e_{\lambda}(t)-1)-1\big)^{k-1}\big|(x)_{n,\lambda}\big\rangle_{\lambda}\nonumber \\
&= \frac{1}{k!}\big\langle e_{\lambda}(t)\big(e_{\lambda}(e_{\lambda}(t)-1)-1\big)^{k-1}\big|(x)_{n,\lambda}\big\rangle_{\lambda}-\frac{1}{k!}\big\langle\big(e_{\lambda}(e_{\lambda}(t)-1)-1\big)^{k-1}\big|(x)_{n,\lambda}\big\rangle_{\lambda}\nonumber\\
&=\frac{1}{k!}\sum_{j=0}^{n-1}\binom{n}{j}(1)_{n-j,\lambda}\sum_{l=0}^{k-1}\binom{k-1}{l}(-1)^{k-1-l}\big\langle e_{\lambda}^{l}(e_{\lambda}(t)-1)\big|(x)_{j,\lambda}\big\rangle_{\lambda} \nonumber\\
&=\frac{1}{k!}\sum_{j=0}^{n-1}\sum_{l=0}^{k-1}\binom{n}{j}\binom{k-1}{l}(1)_{n-j,\lambda}(-1)^{k-1-l}\sum_{i=0}^{\infty}\phi_{i,\lambda}(l)\frac{1}{i!}\big\langle t^{i}|(x)_{j,\lambda}\big\rangle_{\lambda}\nonumber \\
&= \frac{1}{k!}\sum_{j=0}^{n-1}\sum_{l=0}^{k-1}\binom{n}{j}\binom{k-1}{l}(1)_{n-j,\lambda}(-1)^{k-1-l}\phi_{j,\lambda}(l).\nonumber
\end{align}
Therefore, by \eqref{40}, \eqref{41} and \eqref{42}, we obtain the following theorem.
\begin{theorem}
For $n\ge 0$, we have
\begin{displaymath}
\phi_{n,\lambda}(x)=\sum_{l=0}^{n}\beta_{l,\lambda}S_{2,\lambda}(n,l)+\sum_{k=1}^{n}\bigg\{\frac{1}{k!}\sum_{j=0}^{n-1}\sum_{l=0}^{k-1}\binom{n}{j}\binom{k-1}{l}(1)_{n-j,\lambda}(-1)^{k-1-l}\phi_{j,\lambda}(l)\bigg\}b_{k,\lambda}(x).
\end{displaymath}	
\end{theorem}
From Theorem 3, we note that
\begin{align}
d_{m,\lambda}(n,x)&\sim\bigg((mt+1)^{-\frac{1}{m}},\ \frac{1}{\lambda}\big((mt+1)^{\frac{\lambda}{m}}-1\big)\bigg)_{\lambda}\label{43} \\
&\sim\bigg(	(mt+1)^{-\frac{1}{m}},\ \frac{1}{m}\log_{\frac{\lambda}{m}}(mt+1)\bigg)_{\lambda}.\nonumber
\end{align}
For $p(x)\in\mathbb{P}_{n}$, let us assume that
\begin{equation}
p(x)	=\sum_{k=0}^{n}C_{k}d_{m,\lambda}(k,x). \label{44}
\end{equation}
Then we note that
\begin{align}
&\bigg\langle (mt+1)^{-\frac{1}{m}}\bigg(\frac{1}{m}\log_{\frac{\lambda}{m}}(1+mt)\bigg)^{k}\bigg|p(x)\bigg\rangle_{\lambda}	\label{45}\\
&=\sum_{l=0}^{n}C_{l} \bigg\langle (mt+1)^{-\frac{1}{m}}\bigg(\frac{1}{m}\log_{\frac{\lambda}{m}}(1+mt)\bigg)^{k}\bigg|d_{m,\lambda}(l,x)\bigg\rangle_{\lambda}\nonumber\\
&=\sum_{l=0}^{n}C_{l}l!\delta_{l,k}=C_{k}k!.\nonumber
\end{align}
By \eqref{45}, we get
\begin{equation}
C_{k}=\frac{1}{k!} \bigg\langle (mt+1)^{-\frac{1}{m}}\bigg(\frac{1}{m}\log_{\frac{\lambda}{m}}(1+mt)\bigg)^{k}\bigg|p(x)\bigg\rangle_{\lambda}\label{46}.
\end{equation}
Therefore, by \eqref{44} and \eqref{46}, we obtain the following theorem.
\begin{theorem}
For $p(x)\in\mathbb{P}_{n}$, we have
\begin{displaymath}
p(x)=\sum_{k=0}^{n}C_{k}d_{m,\lambda}(k,x),
\end{displaymath}	
where
\begin{displaymath}
C_{k}=\frac{1}{k!} \bigg\langle (mt+1)^{-\frac{1}{m}}\bigg(\frac{1}{m}\log_{\frac{\lambda}{m}}(1+mt)\bigg)^{k}\bigg|p(x)\bigg\rangle_{\lambda}.
\end{displaymath}
\end{theorem}
Let us take $p(x)=\beta_{n,\lambda}(x)$. Then we have
\begin{equation}
\beta_{n,\lambda}(x)=\sum_{k=0}^{n}C_{k}d_{m,\lambda}(k,x),\label{47}	
\end{equation}
where
\begin{align}
C_{k}&=\frac{1}{k!} \bigg\langle (mt+1)^{-\frac{1}{m}}\bigg(\frac{1}{m}\log_{\frac{\lambda}{m}}(1+mt)\bigg)^{k}\bigg|\beta_{n,\lambda}(x)\bigg\rangle_{\lambda}\label{48}\\
&=\sum_{l=0}^{n}\binom{n}{l}\beta_{n-l,\lambda}\frac{1}{k!} \bigg\langle (mt+1)^{-\frac{1}{m}}\bigg(\frac{1}{m}\log_{\frac{\lambda}{m}}(1+mt)\bigg)^{k}\bigg|(x)_{l,\lambda}\bigg\rangle_{\lambda}\nonumber \\
&=\sum_{l=0}^{n}\binom{n}{l}\beta_{n-l,\lambda}\sum_{j=k}^{l}S_{1,\frac{\lambda}{m}}(j,k)m^{j-k}\frac{1}{j!}\langle (mt+1)^{-\frac{1}{m}}t^{j}|(x)_{l,\lambda}\rangle_{\lambda}\nonumber\\
&=\sum_{l=0}^{n}\binom{n}{l}\beta_{n-l,\lambda}\sum_{j=k}^{l}S_{1,\frac{\lambda}{m}}(j,k)\binom{l}{j}m^{j-k}\langle (mt+1)^{-\frac{1}{m}}|(x)_{l-j,\lambda}\rangle_{\lambda}\nonumber \\
&=\sum_{l=0}^{n}\sum_{j=k}^{l}\binom{n}{l}
\beta_{n-l,\lambda}S_{1,\frac{\lambda}{m}}(j,k)m^{j-k}\binom{l}{j}\sum_{i=0}^{\infty}(-1)^{i}\bigg(\frac{1}{m}\bigg)^{i}\frac{1}{i!}\big\langle(\log(mt+1)\big)^{i}|(x)_{l-j,\lambda}\big\rangle_{\lambda}\nonumber\\
&=\sum_{l=0}^{n}\sum_{j=k}^{l}\sum_{i=0}^{l-j}\binom{n}{l}\binom{l}{j}\beta_{n-l,\lambda}S_{1,\frac{\lambda}{m}}(j,k)m^{j-k}(-1)^{i}\bigg(\frac{1}{m}\bigg)^{i}S_{1}(l-j,i)m^{l-j}\nonumber \\
&=\sum_{l=0}^{n}\sum_{j=k}^{l}\sum_{i=0}^{l-j}\binom{n}{l}\binom{l}{j}\beta_{n-l,\lambda}S_{1,\frac{\lambda}{m}}(j,k)m^{l-k-i}(-1)^{i}S_{1}(l-j,i).\nonumber
\end{align}
Therefore, by \eqref{47} and \eqref{48}, we obtain the following theorem.
\begin{theorem}
For $n\ge 0$, we have
\begin{displaymath}
\beta_{n,\lambda}(x)=\sum_{k=0}^{n}\bigg\{\sum_{l=0}^{n}\sum_{j=k}^{l}\sum_{i=0}^{l-j}\binom{n}{l}\binom{l}{j}\beta_{n-l,\lambda}S_{1,\frac{\lambda}{m}}(j,k)m^{l-k-i}(-1)^{i}S_{1}(l-j,i)\bigg\}d_{m,\lambda}(k,x).
\end{displaymath}
\end{theorem}
Let us take $p(x)=(x)_{n,\lambda}\in\mathbb{P}_{n}$. Then, by Theorem 9, we get
\begin{equation}
(x)_{n,\lambda}=\sum_{k=0}^{n}C_{k}d_{m,\lambda}(k,x),\label{49}
\end{equation}
where
\begin{align}
C_{k}&=\frac{1}{k!}\bigg\langle (mt+1)^{-\frac{1}{m}}\bigg(\frac{1}{m}\log_{\frac{\lambda}{m}}(1+mt)\bigg)^{k}\bigg|(x)_{n,\lambda}\bigg\rangle_{\lambda}\label{50}\\
&=\frac{1}{m^{k}}\sum_{l=k}^{n}S_{1,\frac{\lambda}{m}}(l,k)\frac{m^{l}}{l!}\big\langle (mt+1)^{-\frac{1}{m}}t^{l}\big|(x)_{n,\lambda}\big\rangle_{\lambda}\nonumber \\
&=\sum_{l=k}^{n}S_{1,\frac{\lambda}{m}}(l,k)m^{l-k}\binom{n}{l}\big\langle e^{-\frac{1}{m}\log(mt+1)}\big|(x)_{n-l,\lambda}\big\rangle_{\lambda}\nonumber \\
&=\sum_{l=k}^{n}S_{1,\frac{\lambda}{m}}(l,k)m^{l-k}\binom{n}{l}\sum_{i=0}^{n-l}\bigg(-\frac{1}{m}\bigg)^{i}\frac{1}{i!}\big\langle \big(\log(mt+1)\big)^{i}\big|(x)_{n-l,\lambda}\big\rangle_{\lambda}\nonumber 
\end{align}
\begin{align*}
&=\sum_{l=k}^{n}S_{1,\frac{\lambda}{m}}(l,k)m^{l-k}\binom{n}{l}\sum_{i=0}^{n-l}(-1)^{i}\bigg(\frac{1}{m}\bigg)^{i}	\sum_{j=i}^{\infty}S_{1}(j,i)m^{j}\frac{1}{j!}\langle t^{j}|(x)_{n-l,\lambda}\rangle_{\lambda}\nonumber \\
&=\sum_{l=k}^{n}S_{1,\frac{\lambda}{m}}(l,k)m^{l-k}\binom{n}{l}\sum_{i=0}^{n-l}(-1)^{i}m^{-i}S_{1}(n-l,i)m^{n-l} \nonumber \\
&=\sum_{l=k}^{n}\sum_{i=0}^{n-l}\binom{n}{l}S_{1,\frac{\lambda}{m}}(l,k)S_{1}(n-l,i)m^{n-k-i}(-1)^{i}.\nonumber
\end{align*}
Therefore, by \eqref{49} and \eqref{50},we obtain the following theorem.
\begin{theorem}
For $n\ge 0$, we have
\begin{displaymath}
(x)_{n,\lambda}=\sum_{k=0}^{n}\bigg\{\sum_{l=k}^{n}\sum_{i=0}^{n-l}\binom{n}{l}S_{1,\frac{\lambda}{m}}(l,k)S_{1}(n-l,i)m^{n-k-i}(-1)^{i}\bigg\}d_{m,\lambda}(k,x).
\end{displaymath}	
\end{theorem}
Recalling that $d_{m,\lambda}(n,x)\sim\big((mt+1)^{-\frac{1}{m}},\ \frac{1}{m}\log_{\frac{\lambda}{m}}(mt+1)\big)_{\lambda}$, $\phi_{n,\lambda}(x)\sim (1,\log_{\lambda}(1+t))_{\lambda}$, \\
let us assume that
\begin{equation}
d_{m,\lambda}(n,x)=\sum_{k=0}^{n}C_{n,k}\phi_{k,\lambda}(x).\label{51}	
\end{equation}
Then, by \eqref{9}, we get
\begin{align}
C_{n,k}&=\frac{1}{k!}\bigg\langle e_{\lambda}(t)\bigg(\log_{\lambda}\bigg(1+\frac{1}{m}(e_{\lambda}^{m}(t)-1)\bigg)\bigg)^{k}\bigg|(x)_{n,\lambda}\bigg\rangle_{\lambda}\label{52}	\\
&= \sum_{j=k}^{n}S_{1,\lambda}(j,k)\frac{1}{j!}\bigg\langle e_{\lambda}(t)\bigg(\frac{1}{m}\big(e_{\lambda}^{m}(t)-1\big)\bigg)^{j}\bigg|(x)_{n,\lambda}\bigg\rangle_{\lambda}\nonumber \\
&=\sum_{j=k}^{n}S_{1,\lambda}(j,k)\sum_{l=j}^{\infty}W_{m,\lambda}(l,j)\frac{1}{l!}\langle t^{l}|(x)_{n,\lambda}\rangle_{\lambda}\nonumber \\
&=\sum_{j=k}^{n}S_{1,\lambda}(j,k)W_{m,\lambda}(n,j).\nonumber
\end{align}
Therefore, by \eqref{51} and \eqref{52}, we obtain the following theorem.
\begin{theorem}
For $n\ge 0$, we have
\begin{displaymath}
d_{m,\lambda}(n,x)=\sum_{k=0}^{n}\bigg\{\sum_{j=k}^{n}S_{1,\lambda}(j,k)W_{m,\lambda}(n,j)\bigg\}\phi_{k,\lambda}(x).
\end{displaymath}	
\end{theorem}
 Now, we observe that
 \begin{align}
 &e_{\lambda}^{x}\bigg(\frac{1}{m}\big(e_{\lambda}^{m}(t)-1\big)\bigg)=\sum_{k=0}^{\infty}(x)_{k,\lambda}\bigg(\frac{1}{m}\bigg)^{k}\frac{1}{k!}\big(e_{\frac{\lambda}{m}}(mt)-1\big)^{k}\label{53} \\
 &=\sum_{k=0}^{\infty}(x)_{k,\lambda}m^{-k}\sum_{n=k}^{\infty}S_{2,\frac{\lambda}{m}}(n,k)m^{n}\frac{t^{n}}{n!}\nonumber \\
 &=\sum_{n=0}^{\infty}\bigg(\sum_{k=0}^{n}\bigg(\frac{x}{m}\bigg)_{k,\frac{\lambda}{m}}S_{2,\frac{\lambda}{m}}(n,k)\bigg)m^{n}\frac{t^{n}}{n!}=\sum_{n=0}^{\infty}\phi_{n,\frac{\lambda}{m}}\bigg(\frac{x}{m}\bigg)m^{n}\frac{t^{n}}{n!}. \nonumber
 \end{align}
From \eqref{7} and \eqref{53}, we note that
\begin{equation}
m^{n}\phi_{n,\frac{\lambda}{m}}\bigg(\frac{x}{m}\bigg)\sim\bigg(1,\frac{(mt+1)^{\frac{\lambda}{m}}-1}{\lambda}\bigg)_{\lambda}\sim \bigg(1,\frac{1}{m}\log_{\frac{\lambda}{m}}(mt+1)\bigg)_{\lambda}.\label{54}	
\end{equation}
Let us assume that
\begin{equation}
m^{n}\phi_{n,\frac{\lambda}{m}}\bigg(\frac{x}{m}\bigg)=\sum_{k=0}^{n}C_{n,k}d_{m,\lambda}(k,x).\label{55}
\end{equation}
Then, by \eqref{9}, \eqref{43} and \eqref{54}, we get
\begin{align}
C_{n,k}=\frac{1}{k!}\bigg\langle e_{\lambda}^{-1}(t) t^{k}\bigg|(x)_{n,\lambda}\bigg\rangle_{\lambda} 
=\binom{n}{k} \bigg\langle e_{\lambda}^{-1}(t) \bigg|(x)_{n-k,\lambda} \bigg \rangle_{\lambda} =\binom{n}{k}(-1)_{n-k,\lambda}.\label{56} 
\end{align}
Thus, by \eqref{55} and \eqref{56}, we get
\begin{displaymath}
\phi_{n,\frac{\lambda}{m}}\bigg(\frac{x}{m}\bigg)=m^{-n}\sum_{k=0}^{n}\binom{n}{k}(-1)_{n-k,\lambda}d_{m,\lambda}(k,x).
\end{displaymath}

\section{Conclusion}
In this paper, we introduced the fully degenerate Bell polynomials $\phi_{n,\lambda}(x)$, which are degenerate versions of the ordinary Bell polynomials $\mathrm{Bel}_{n}(x).$ Here  $\mathrm{Bel}_{n}(x)$ is the natural extension of the Bell number $\mathrm{Bel}_{n}$, which
counts the number of partitions of a set with $n$ elements into disjoint nonempty subsets. Further, we introduced the fully degenerate Dowling polynomials $d_{m,\lambda}(n,x)$ which are degenerate versions of the ordinary Dowling polynomials $D_{m}(n,x)$. Here $D_{m}(n,x)$ is the natural extension of the Whitney numbers of the second kind $W_{m}(n,k)$, which is equal to the number of coloured partitions of $\left\{1, \dots, n + 1 \right\}$ into $k + 1$ nonempty subsets satisfying the conditions stated in Introduction with $r=1$. \par
We derived some identities and properties for those degenerate polynomials by using the recently introduced $\lambda$-umbral calculus. Especially, we derived formulas expressing any polynomial in terms of the fully degenerate Bell and the fully degenerate Dowling polynomials and applied 
them to the degenerate Bernoulli polynomials and degenerate falling factorial polynomials. \par
It is one of our future projects to continue to work on various degenerate versions of many special numbers and polynomials.

\end{document}